\documentclass[11pt]{article}
\usepackage{a4wide,parskip,amsmath,amssymb,latexsym,graphicx}
\usepackage{amsmath,amsfonts,amssymb,latexsym,graphics,epsfig,url}
\usepackage{color}
\usepackage{amsthm,enumerate}
\usepackage[english]{babel}
\usepackage{graphicx}
\usepackage{bm}
\usepackage{xcolor}

\newtheorem{theorem}{Theorem}
\newtheorem{proposition}{Proposition}
\newtheorem{lemma}{Lemma}
\newtheorem{definition}{Definition}
\newtheorem{example}{Example}
\newtheorem{corollary}{Corollary}

\def\M{{\bm{M}}}
\def\S{\mbox{\boldmath $S$}}
\def\U{\mbox{\boldmath $U$}}
\def\V{\mbox{\boldmath $V$}}
\def\o{\bm{0}}
\def\m{\bm{m}}
\def\llambda{\bm{\lambda}}
\newcommand{\sq}[1]{[\![{#1}]\!]}
\newcommand{\sg}[1]{\langle{#1}\rangle}
\newcommand{\norm}[1]{\|#1\|_1}
\newcommand{\mmod}[1]{\!\!\pmod{#1}}
\def\LL{\mathcal{L}}
\def\HH{\mathcal{H}}

\def\N{\mathbb{N}}
\def\Z{\mathbb{Z}}
\def\R{\mathbb{R}}
\def\b{\mbox{\boldmath $b$}}
\def\e{\mbox{\boldmath $e$}}
\def\c{\mbox{\boldmath $c$}}
\def\u{{\mbox {\boldmath $u$}}}
\def\v{{\mbox {\boldmath $v$}}}

\def\x{\mbox{\boldmath $x$}}
\def\y{\mbox{\boldmath $y$}}

\def\a{\mbox{\boldmath $a$}}
\def\b{\mbox{\boldmath $b$}}

\def\Cay{\mathop{\rm Cay }\nolimits}

\def\diag{\mathop{\rm diag }\nolimits}

\newcommand{\NA}{\textrm{NA}}
\newcommand{\NC}{\textrm{NC}}
\newcommand{\lb}{\textrm{lb}}

\begin{document}

\title{The Dilating Method to obtain dense \\
Cayley digraphs on finite Abelian groups
\thanks{
	This research has been partially supported by
	AGAUR from the Catalan Government under project 2017SGR1087, and by MICINN from the Spanish Government with the European Regional Development Fund under project
	PGC2018-095471-B-I00.
\newline \indent
 E-mails: {\tt francesc.aguilo@upc.edu, miguel.angel.fiol@upc.edu, sonia.perez-mansilla@upc.edu}
            }
}
\author{F. Aguil\'o$^a$, M. A. Fiol$^b$ and S. P\'erez$^c$\\
\\ {\small $^{a,b,c}$Departament de Matem\`atiques}
\\ {\small Universitat Polit\`ecnica de Catalunya}
\\ {\small Jordi Girona 1-3 , M\`odul C3, Campus Nord }
\\ {\small 08034 Barcelona, Catalonia (Spain).}
\\ {\small $^b$Barcelona Graduate School of Mathematics}
\\ {\small Institut de Matem\`atiques de la UPC-BarcelonaTech (IMTech)}
\\ {\small Barcelona, Catalonia.}
}

\maketitle
\begin{abstract}
A geometric method for obtaining an infinite family of Cayley digraphs of constant density on finite Abelian groups is presented. The method works for any given degree $d\ge 2$, and it is based on consecutive dilates of a minimum distance diagram associated with a given initial Cayley digraph.
The method is used to obtain infinite families of dense or asymptotically dense Cayley digraphs. In particular, for degree $d=3$, the first explicit infinite family of maximum known density $\delta=0.084$ is obtained.
\end{abstract}



\section{Introduction}
Some discrete problems are studied by many authors using bare geometrical methods. Some of them are worth mentioning here because of the common used tool. Namely, the so-called {\em minimum distance diagrams} (MDD's), which are defined in the next section. Some examples are in the study of {\em numerical semigroups}, that is the {\em Frobenius number}, the {\em set of factorizations}, and the {\em denumerant} \cite{AgGa:10,BeNiHeWa:12}. Also, some metric properties of {\em Cayley digraphs} of finite Abelian groups  have been studied using MDD's, mainly the {\em diameter} and the {\em density} \cite{AgFiGa:97,DoFa:04,FiFoZi:98,Rod96}. Cayley digraphs of finite Abelian groups will be studied in this work.

MDD's were introduced by Wong and Coppersmith in 1974 \cite{wc74} for degree $d=2$. Fiol, Yebra, Alegre, and Valero in 1987 \cite{fyav87} gave infinite families of digraphs of degree $d=2$ with optimal diameter. Their results were based on geometric arguments of the related MDD's. These results were generic for degree $d=2$. The Smith normal form was used to obtain the related generators. A digraph of degree $d=3$ with optimal diameter $k=8$ was also given. This is the Cayley digraph on the cyclic group  $\Z_{111}$ and generating set $\{1,31,69\}$. It was found by computer search, and a related MDD was also depicted. Infinitely many infinite families of optimal diameter digraphs of degree $d=2$ were given in \cite{eaf93} using the same geometrical approach.

Aguil\'o, Fiol, and Garc\'{\i}a in 1997 \cite{AgFiGa:97} gave some infinite families of Cayley digraphs of degree $d=3$ using three-dimensional MDD's and the Smith normal form. The underlying idea for obtaining an infinite family of digraphs with good metrical properties is sketched here. From an initial digraph with good diameter, take the related MDD $\HH_0$. Then, the `dilation' of the cubes that belong to $\HH_0$ gives another MDD $\HH_1$ whose sides preserve their relative ratios. This approach seems to give a good diameter to the related digraph. All the details when looking for the generators, the resulting diameter, and related parameters depended on the particular initial $\HH_0$. However, neither a generic method for degree $d=3$, nor for degree $d\geq 4$ was provided.

Fiduccia, Forcade and Zito in 1998 \cite{FiFoZi:98} gave nice results on diameter and density for degree $d=3$. In their work, we can find many geometrical results on MDD's. In particular, these authors use the digraph diameter $k$ and the `{\em solid diameter}' $D$ of the related MDD ($D=k+3$). Perhaps this is the first work where geometrical generic results on three-dimensional MDD were given. In particular, an upper bound for the three dimensional MDD's density, defined as $\delta=\textrm{Vol}/D^3$, is also given ($\delta\leq 3/25 =0.12$) in \cite[Sec.~1]{FiFoZi:98}. The authors pointed out that the best density $\delta_0=0.084$ for an MDD $\HH_1$ of volume $n_1=84$ is the same as the best density $\delta_2=\delta_0$ of another MDD $\HH_2$ of volume $n_2=672=84\times2^3$. The tile $\HH_2$ can be seen as a dilation of $\HH_1$ (each unitary cube of $\HH_1$ is dilated by a $2\times2\times2$ cube). The same is true for $\HH_3$ and $\HH_4$ of volumes $n_3=2268=84\times3^3$ and $n_4=5376=84\times4^3$, respectively (similar dilatations $3\times3\times3$ and $4\times4\times4$, respectively). The idea of `dilating' the unit cubes of a given MDD to obtain an infinite family of MDD's of the same density is commented for $d=3$ in \cite[Sec.~5]{FiFoZi:98}. No proof is given, however. Besides, the resulting generators for the related Cayley digraphs are not mentioned. In the `{\em Open problems}' section of the same work, several questions were posed. In particular, it was asked whether the results could be generalized to four or more dimensions. Forcade and Lamoreaux in 2000 \cite{FoLa:00} proved that the maximum density for degree $d=2$ is $\Delta_2=\frac13$ and also, when $d=3$, the density has a local maximum at $\Delta_3'=\delta_0=0.084$. As far as we know, no Cayley digraph of a finite Abelian group with density larger than $\delta_0$ has been found.

Another excellent paper on this topic is the work of Dougherty and Faber in 2004 \cite{DoFa:04}. They also use a geometrical approach by associating tiles to the directed and undirected Cayley graphs. A table of several optimal digraphs, found by computer search, is included for degree $d=3$. In particular, there are some digraphs of density $\delta_0$.
Several asymptotic results are also given for degree $d=3$ and higher degrees.
However, from the results in \cite[pag. 502]{DoFa:04}, no generic method for the dilating process in the case of degree $d=3$ or larger is discussed.

In this work we introduce the {\em Dilating Method} based on minimum distance diagrams. It works for any degree $d\geq2$ and gives the generators that allow finding explicit constructions of Cayley digraphs with constant density.  We use it to derive a dense infinite family of Cayley digraphs of degree $d=3$ and constant density $\delta_0=0.084$. As far as we know, this is the first infinite family of tiles (or digraphs) with this density. An infinite set of asymptotically dense families of  Cayley digraphs is also explicitly built for any degree $d\geq3$.

An overview of the contents of the paper is as follows: Section~2 contains the main definitions and preliminary results. The Dilating Method, together with its proof, is presented in Section~3. As commented above, this method is valid for any degree $d\geq 2$. As far as we know, it is the first generic geometrical method in this context. Some new dense and asymptotically dense families are constructed in the last section as an application of the method. More precisely, when applied to the degree $d=3$, an {\bf explicit} infinite family of digraphs of density $\delta=0.084$ is found. To our knowledge, until now only a few digraphs having this density were known.
Concerning the employed examples, they are the same (or closely related to) ones already used in other papers. See, for instance, \cite{AgSiZa:01}. The reason is that they constitute known ``optimal cases" and so, they appear when an exhaustive computer search is carried out.

\section{Preliminaries}
\label{preli}
Consider an integral matrix $\M\in\Z^{n\times n}$ with $N=|\det\M|\neq 0$, and with \textit{Smith normal form} $\S=\U\M\V$, for unimodular matrices $\U,\V\in\Z^{n\times n}$ (see, for example, \cite{Ne:72}). Let $\M\Z^n$ denote the lattice generated by the columns of $\M$ which, with the usual vector addition, is a normal subgroup of $\Z^n$. Let $\Z^n/\M\Z^n$ be the quotient group induced by
the equivalence relation in $\Z^n$ given by
$$
\a\equiv\b\; \mmod{\M}\quad\Leftrightarrow\quad\a-\b\in \M\Z^n \quad \mbox{(that is, $\exists\llambda\in\Z^n:\;\a-\b=\M\llambda$)}.
$$
Let us consider the Cayley digraph $G_{\M}=\Cay(\Z^n/\M\Z^n,E_n)$ where $E_n=\{\e_1,\ldots,\e_n\}$ is the canonical basis of $\Z^n$.
Let $[r,s)=\{x\in\R:r\leq x<s\}$. Given $\a=(a_1,\ldots,a_n)\in\Z^n$, we denote the unitary cube $\sq{\a}=\sq{a_1,\ldots,a_n}=[a_1,a_1+1)\times\cdots\times[a_n,a_n+1)\subset\R^n$. In this sense, the cube $\sq{\a}$ represents the vertex $(a_1,\ldots,a_n)$ in $G_\M$ with the equivalence relation
$\sq{\a}\sim\sq{\b}$ whenever $\a\equiv\b\; \mmod{\M}$ and $\sq{\a}\not\sim\sq{\b}$ otherwise.
For more details about congruences in $\Z^n$ and their role in the study of Cayley digraphs on Abelian groups, see \cite{eaf93,Fi:87,Fi:95}.

For a given pair $\a,\b\in\Z^n$, we denote by $\a\leq\b$ when the inequality holds for each coordinate. Let $\N=\Z_{\geq0}$ denote the set of nonnegative integers. Given $\a\in\N^n$, consider the set of unitary cubes $\nabla(\a)=\{\sq{\b}:\o\leq\b\leq\a\}$.
\\
\begin{definition}[Hyper-L]
\label{def:hyperL}
Given a finite Abelian group $\Gamma=\sg{\gamma_1,\ldots,\gamma_n}$ of order $N=|\Gamma|$, consider the map $\phi:\N^n\longrightarrow\Gamma$ given by $\phi(\a)=a_1\gamma_1+\cdots+a_n\gamma_n$. A {\em hyper-L} of the Cayley digraph $\Cay(\Gamma,\{\gamma_1,\ldots,\gamma_n\})$, denoted by $\LL$, is a set of $N$ unitary cubes $\LL=\{\sq{\a_0},\ldots,\sq{\a_{N-1}}\}$ such that
\begin{itemize}
\item[$(i)$] $\{\phi(\a):\sq{\a}\in\LL\}=\Gamma$,
\item[$(ii)$] $\sq{\a}\in\LL\Rightarrow\nabla(\a)\subset\LL$.
\end{itemize}
\end{definition}

Given $\x\in\Z^n$, let us consider the $\ell_1$ norm $\norm{\x}=|x_1|+\cdots+|x_n|$. The \textit{diameter of the hyper-L} $\LL$ is defined to be $k(\LL)=\max\{\norm{\a}:\sq{\a}\in\LL\}$. For the usual definition of the diameter  $k(G)$ of a Cayley digraph $G=\Cay(\Gamma,\sg{\gamma_1,\ldots,\gamma_n})$, we have $k(G)\leq k(\LL)$ for every hyper-L $\LL$ of $G$.
\\
\begin{definition}[Minimum distance diagram]
\label{def:MDD}
A {\em minimum distance diagram} $\HH$ of the Cayley digraph $G$ is a hyper-L satisfying $\norm{\a}=\min\{\norm{\x}:\x\in\phi^{-1}(\phi(\a))\}$ for all $\sq{\a}\in\HH$.
\end{definition}

It follows that $k(G)=k(\HH)$ for each minimum distance diagram $\HH$ of $G$, and also $k(G)=\min\{k(\LL):\textrm{$\LL$ is a hyper-L of $G$}\}$. When $\Gamma$ is a cyclic group, Definition~\ref{def:MDD} is equivalent to \cite[Def.~2.1]{SaSa:09} for multiloop networks.

\begin{figure}[h]
\centering
\includegraphics[width=0.45\linewidth]{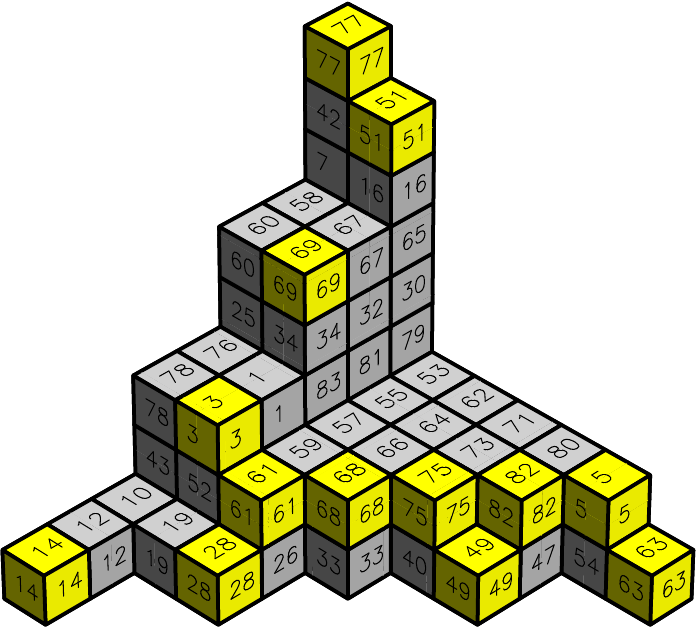}
\hspace{0.07\linewidth}
\includegraphics[width=0.42\linewidth]{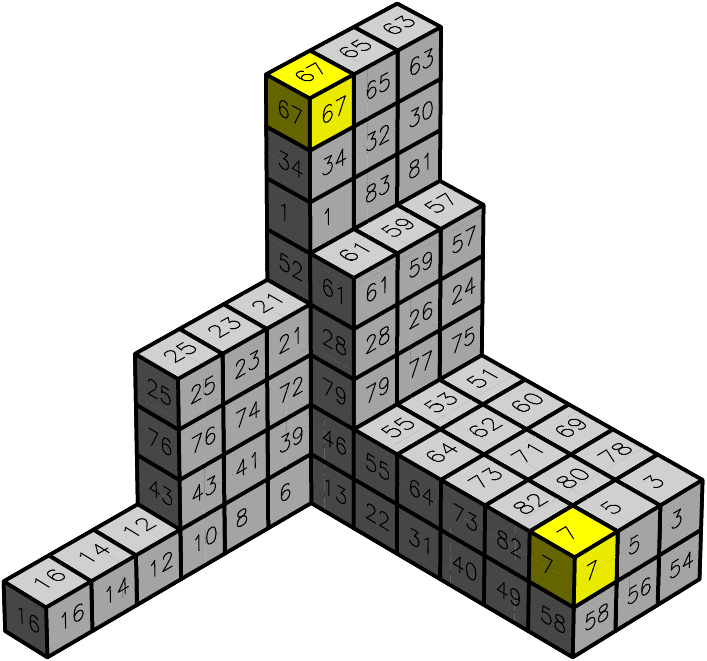}
\caption{MDD's of $G_1=\Cay(\Z_{84},\{2,9,35\})$ and $G_2=\Cay(\Z_{84},\{2,9,33\})$.}
\label{fig:L84}
\end{figure}

\begin{example}
\label{ex:L84}
The Cayley digraph $G_1=\Cay(\Z_{84},\{2,9,35\})$ has only one minimum distance diagram $\HH_1$ which is shown on the left hand side of Figure~\ref{fig:L84}. The diameter is $k(G_1)=k(\HH_1)=\norm{(0,0,7)}=7$. The cube $\sq{0,0,7}$ corresponds to vertex $77$ in $G_1$. There are 13 unit cubes in $\HH_1$ with maximum norm. Namely, $\sq{0,0,7}$, $\sq{0,1,6}$, $\sq{2,1,4}$, $\sq{4,1,2}$, $\sq{0,6,1}$, $\sq{1,5,1}$, $\sq{2,4,1}$, $\sq{3,3,1}$, $\sq{4,2,1}$, $\sq{0,7,0}$,  $\sq{2,5,0}$, $\sq{5,2,0}$, and $\sq{7,0,0}$. These cubes are painted in light color. Each cube $\a=\sq{a_1,a_2,a_3}$ in $\HH_1$ corresponds to the vertex $2a_1+9a_2+35a_3$ in $\Z_{84}$.

The minimum distance diagram $\HH_2$ shown on the right hand side of Figure~\ref{fig:L84} has diameter $k(\HH_2)=9$, and corresponds to the digraph $G_2=\Cay(\Z_{84},\{2,9,33\})$. It has only two cubes with maximum norm: $\sq{2,6,1}$ and $\sq{2,0,7}$. They correspond to the vertices $7$ and $67$, respectively. Exhaustive computer search shows that the minimum diameter for a Cayley digraph on a three generated cyclic group with order $84$ is $7$.
\end{example}

For a Cayley digraph $G$ of an Abelian group, with order $N(G)$, degree $d(G)$, and diameter $k(G)$, Fiduccia, Forcade and Zito  \cite{FiFoZi:98} defined its \textit{solid density} as
\begin{equation}
\label{density}
\delta(G)=\frac{N(G)}{(k(G)+d(G))^{d(G)}},
\end{equation}
the density of an MDD related to $G$. From now on, we use this expression for computing the density. Let us denote
\begin{equation*}
\Delta_{d,k}=\max\{\delta(G): d(G)=d, k(G)=k\},\quad\mbox{and}\quad
\Delta_d=\max\{\Delta_{d,k}: d(G)=d\}.
\end{equation*}
Forcade and Lamoreaux proved that $\Delta_2=1/3$ by using a generic optimal diagram $\HH$ in \cite[Section~4]{FoLa:00}. This diagram was already known by Fiol, Yebra, Alegre, and Valero \cite{fyav87} to derive the tight lower bound $k(N)\geq\lb(N)=\lceil\sqrt{3N}\rceil-2$ for a generic Cayley digraph of degree $d=2$ and order $N$. The density $\Delta_2$ is attained by the digraphs $G_{\M}=\Cay(\Z^2/\M\Z^2,E_2)$, where $\M$ is the matrix with rows $(2t,-t)$ and $(-t,2t)$ or, equivalently, $G_t=\Cay(\Z_t\oplus\Z_{3t},\{(1,-1),(0,1)\})$, with $N_t=3t^2$  and $k(G_t)=3t-2$ for each $t\geq1$ (see next section). No Cayley digraph on cyclic group $\Z_{3t^2}$ of degree $d=2$ attains this density for $t>1$.

As far as we know, $\Delta_d$ remains unknown for $d>2$. For $d=3$ there are some main facts:
Fiduccia, Forcade, and Zito in \cite[Cor.~3.6]{FiFoZi:98} proved that $\Delta_3\leq3/25=0.12$. However, numerical evidences  suggest that $\Delta_3$ would be a smaller value. The maximum density attained by known Cayley digraphs is $\delta_0=0.084$ and they have been found by computer search. Such digraphs are $F_0$, $F_1$ in \cite[Table~1]{FiFoZi:98}, and $F_1'\cong F_0$, $F_2$, $F_3$ in \cite[Table~8.2]{DoFa:04}, where:
\begin{itemize}
\item
$F_0\cong G_1=\Cay(\Z_{84},\{2,9,35\})$;
\item
$F_1=\Cay(\Z_2\oplus\Z_2\oplus\Z_{168},\{(1,0,2),(0,0,9),(0,1,35)\})$;
\item
$F_2=\Cay(\Z_3\oplus\Z_3\oplus\Z_{252},\{(0,0,2),(0,1,9),(1,0,35)\})$;
\item
$F_3=\Cay(\Z_4\oplus\Z_4\oplus\Z_{336},\{(0,1,2),(0,0,9),(1,0,35)\})$.
\end{itemize}
Let $\NA_{d,k}$ (respectively, $\NC_{d,k}$) be the maximum number of vertices that a Cayley digraph of an Abelian group (respectively, of a cyclic group), with  degree $d$ and diameter $k$, can have. Let us denote $\lb(d,k)$ the {\em lower bound} for $\NA_{d,k}$. Then, from Wong and Coppersmith \cite{wc74} and  Dougherty and Faber \cite[Theorem~9.1]{DoFa:04}, it follows that, for $d>1$,
\begin{equation}
\lb(d,k)=\frac{c}{d(\ln d)^{1+\log_2 e}}\frac{k^d}{d!}+O(k^{d-1})\leq\NA_{d,k}<{k+d\choose k},\label{eq:df}
\end{equation}
for some  constant $c$. As far as we know, no constructions of Cayley digraphs $G$ of order $N(G)\sim\lb(d,k)$ are known.

Notice that, for a given Cayley digraph $G$, a large value of the ratio $N(G)/k(G)$ does not guarantee a high density of $G$. In this work we propose the {\em Dilating Method} which allows the generation of an infinite family of dense digraphs from a given initial dense digraph.

\section{The Dilating Method}
For a given finite Abelian Cayley digraph $G=\Cay(\Gamma,\{\gamma_1,\ldots,\gamma_n\})$ of degree $d=n$, with minimum distance diagram $\HH$, there is an integral matrix $\M\in\Z^{n\times n}$ such that
\begin{equation}
G\cong\Cay(\Z^n/\M\Z^n,E_n)\cong\Cay(\Z_{s_1}\oplus\cdots\oplus\Z_{s_n},\{\u_1,\ldots,\u_n\}),\label{eq:isom}
\end{equation}
where $\S=\diag(s_1,\ldots,s_n)=\U\M\V$ is the \textit{Smith normal form} of $\M$, $E_n=\{\e_1,\ldots,\e_n\}$ is the canonical basis of $\Z^n$ and $\{\u_1,\ldots,\u_n\}$ are the column vectors of the matrix $\U$. Let $\HH$ be a minimum distance diagram of $G$. It is well known that $\HH$ tessellates $\R^n$ by translation through the vector set $C_\M=\{\m_1,\ldots,\m_n\}$, with $C_{\M}$ being the  set of column vectors of $\M$. The map $\psi(\x)=\U\x$ plays an important role because of the equivalence $$
\a\equiv\b\;\mmod{\M}\quad \Leftrightarrow\quad \psi(\a)\equiv\psi(\b)\;\mmod{\S}
$$ and so, it follows that  $\sq{\a}\sim\sq{\b}\Leftrightarrow\psi(\a)\equiv\psi(\b)$ in $\Z_{s_1}\oplus\cdots\oplus\Z_{s_n}$. For more details, see \cite{eaf93,Fi:95}.
\\
\begin{example}
\label{ex:M84}
Let us consider again the Cayley digraph $G_1$ and its minimum distance diagram $\HH_1$ of Example~\ref{ex:L84}. In this case we have $\S=\diag(1,1,84)$ and
\[
\U\M\V=
\left(\begin{array}{rrr}
0 & 0 & 1 \\
-2 & 1 & 10 \\
7 & -3 & -38
\end{array}\right)
\left(\begin{array}{rrr}
1 & 2 & -6 \\
5 & 2 & 4 \\
2 & -2 & 3
\end{array}\right)
\left(\begin{array}{rrr}
2 & 1 & 26 \\
0 & 1 & 23 \\
-1 & 0 & -2
\end{array}\right)=\S.
\]
Thus, $\HH_1$ tessellates $\R^3$ through the lattice generated by $\{(1,5,2)^\top,(2,2,-2)^\top,(-6,4,3)^\top\}$ and $G_1\cong\Cay(\Z^3/\M\Z^3,E_3)$. It is not difficult to see by computer that $\psi(\a)\not\equiv\psi(\b)$ in $\Z_1\oplus\Z_1\oplus\Z_{84}\cong\Z_{84}$ for any pair of cubes $\sq{\a},\sq{\b}\in\HH_1$. The required optimality condition of Definition~\ref{def:MDD} can also be verified on $\HH_1$ by computer. Notice that $G_1\cong\Cay(\Z_{84},\{46,81,7\})$ by the isomorphism $\mu(x)=53x$ since $46\equiv-38\;\mmod{84}$, $81\equiv-3\;\mmod{84}$, and $\gcd(84,53)=1$.
The choice of the optimal lattice is not unique. For instance, Daugherty and Faber  \cite{DoFa:04} used the generating set
$\{\v_1,\v_2,\v_3\}=\{(-2,2,2)^\top,(3,-3,3)^\top,(4,3,-1)^\top\}$.  Then, taking into account that the following $($column$)$ vectors belong also to the lattice:
{\small
\begin{align*}
\v_4 &=(6,1,-3)=\v_3-\v_1,\quad \v_5=(5,-5,1)=\v_2-\v_1,\quad \v_6=(1,6,-4)=\v_3-\v_2,\\
\v_7 &=(2,5,1)=\v_1+\v_3,\quad \v_8=(1,-1,5)=\v_1+\v_2,\quad \v_9=(-1,1,7)=2\v_1+\v_2,\\
\v_{10} &=(3,4,-6)=\v_3-\v_1-\v_2,\quad \v_{11}=(-7,7,1)=2\v_1-\v_2,\quad \v_{12}=(-5,-2,8)=2\v_1+\v_2-\v_3,\\
\v_{13} &=(-1,8,-2)=\v_1-\v_2+\v_3,\quad \v_{14}=(8,-1,-5)=\v_3-2\v_1,
\end{align*}
}
it is not difficult to check by hand that the obtained hyper-$L$ represents a digraph with  diameter 7. For instance, following a similar approach than in \cite{DoFa:04},  Figure \ref{farthest-vertices} shows that all integral (column) vectors $\x=x_1x_2x_3=(x_1,x_2,x_3)$ at distance eight  from the origin $000$ (that is, all integral points satisfying $x_1+x_2+x_3=8$) are equivalent, modulo $\M=(\v_1|\v_2|\v_3)$, with vectors which are at distance $\le 7$. In particular, notice that $(2,5,1)\equiv (0,0,0)$.
\begin{figure}[h!]
{\footnotesize
\begin{center}
$\begin{array}{c}
008-\v_{12} \\
= 520
\end{array}$

$\begin{array}{c}
107-\v_{8} \\
= 012
\end{array}$
$\begin{array}{c}
017-\v_{9} \\
= 100
\end{array}$

$\begin{array}{c}
206-\v_{8} \\
= 111
\end{array}$
$\begin{array}{c}
116-\v_{8} \\
= 021
\end{array}$
$\begin{array}{c}
026-\v_{1} \\
= 204
\end{array}$

$\begin{array}{c}
305-\v_{2} \\
= 032
\end{array}$
$\begin{array}{c}
215-\v_{8} \\
= 120
\end{array}$
$\begin{array}{c}
125-\v_{1} \\
= 303
\end{array}$
$\begin{array}{c}
035-\v_{1} \\
= 213
\end{array}$

$\begin{array}{c}
404-\v_{2} \\
= 131
\end{array}$
$\begin{array}{c}
314-\v_{2} \\
= 041
\end{array}$
$\begin{array}{c}
224-\v_{1} \\
= 402
\end{array}$
$\begin{array}{c}
134-\v_{1} \\
= 312
\end{array}$
$\begin{array}{c}
044-\v_{1} \\
= 222
\end{array}$

$\begin{array}{c}
503-\v_{2} \\
= 230
\end{array}$
$\begin{array}{c}
413-\v_{2} \\
= 140
\end{array}$
$\begin{array}{c}
323-\v_{1} \\
= 501
\end{array}$
$\begin{array}{c}
233-\v_{1} \\
= 411
\end{array}$
$\begin{array}{c}
143-\v_{1} \\
= 321
\end{array}$
$\begin{array}{c}
053-\v_{1} \\
= 231
\end{array}$

$\begin{array}{c}
602-\v_{5} \\
= 151
\end{array}$
$\begin{array}{c}
512-\v_{5} \\
= 061
\end{array}$
$\begin{array}{c}
422-\v_{1} \\
= 600
\end{array}$
$\begin{array}{c}
332-\v_{1} \\
= 510
\end{array}$
$\begin{array}{c}
242-\v_{1} \\
= 420
\end{array}$
$\begin{array}{c}
152-\v_{1} \\
= 330
\end{array}$
$\begin{array}{c}
062-\v_{1} \\
= 240
\end{array}$

$\begin{array}{c}
701-\v_{5} \\
= 250
\end{array}$
$\begin{array}{c}
611-\v_{4} \\
= 004
\end{array}$
$\begin{array}{c}
521-\v_{5} \\
= 070
\end{array}$
$\begin{array}{c}
431-\v_{3} \\
= 002
\end{array}$
$\begin{array}{c}
341-\v_{10} \\
= 007
\end{array}$
$\begin{array}{c}
251-\v_{7} \\
= {\bf 000}
\end{array}$
$\begin{array}{c}
161-\v_{6} \\
= 005
\end{array}$
$\begin{array}{c}
071-\v_{11} \\
=700
\end{array}$

$\begin{array}{c}
800-\v_{14} \\
= 015
\end{array}$
$\begin{array}{c}
710-\v_{4} \\
= 103
\end{array}$
$\begin{array}{c}
620-\v_{4} \\
= 013
\end{array}$
$\begin{array}{c}
530-\v_{3} \\
= 101
\end{array}$
$\begin{array}{c}
440-\v_{3} \\
= 011
\end{array}$
$\begin{array}{c}
350-\v_{10} \\
= 016
\end{array}$
$\begin{array}{c}
260-\v_{6} \\
= 104
\end{array}$
$\begin{array}{c}
170-\v_{6} \\
= 014
\end{array}$
$\begin{array}{c}
080-\v_{13} \\
= 102
\end{array}$
\end{center}
}
\caption{Every  integer vector $\x=(x_1,x_2,x_3)=x_1x_2x_3$  at distance 8 from the origin is congruent modulo $\M=(\v_1|\v_2|\v_3)$ with $\y=y_1y_2y_3$ at distance $\le 7$. The points $\x$ are in the intersection of the plane $x_1+x_2+x_3=8$ with the planes limiting the first octant.}
\label{farthest-vertices}
\end{figure}
\end{example}

Given a unit cube $\sq{\a}$, consider the dilate $t\sq{\a}$ defined by
\begin{equation}
t\sq{\a}=\{\sq{t\a+(\alpha_1,\ldots,\alpha_n)}:0\leq \alpha_1,\ldots,\alpha_n\leq t-1\}\label{eq:cubdil}
\end{equation}
for $t\geq1$, $t\in\N$. Notice that $t\sq{\a}\cap t\sq{\b}=\emptyset$ whenever $\sq{\a}\neq\sq{\b}$.
\\
\begin{definition}[The dilation of a hyper-L]
\label{def:dilate}
Let $\LL$ be a hyper-L of some Cayley digraph. The $t$--dilate $t\LL$ of $\LL$ is defined by $t\LL=\{t\sq{\a}:\sq{\a}\in\LL\}$.
\end{definition}

This definition corresponds to a dilatation of $\LL$ in such a way that each unit cube $\sq{\a}$ in $\LL$ is dilated to another cube, $t\sq{\a}$, of side $t$ in $t\LL$ according to \eqref{eq:cubdil}.

\begin{lemma}
\label{lem:grid}
Given a hyper-L $\LL$ of the digraph $G_{\M}$, consider the dilate $t\LL$ for a given integer $t\geq1$. Then, $\x\not\equiv\y\;\mmod{t\M}$ for any pair of cubes $\sq{\x},\sq{\y}\in t\LL$.
\end{lemma}

\begin{proof}
Assume that $\sq{\x}$ and $\sq{\y}$ are different unit cubes in $t\LL$, which can be written in the form   $\x=t\a+(\alpha_1,\ldots,\alpha_n)^\top$ and $\y=t\b+(\beta_1,\ldots,\beta_n)^\top$, with  $\sq{\a},\sq{\b}\in\LL$,  $0\leq\alpha_i,\beta_i<t$ for $1\leq i\leq n$, and $\alpha_1+\cdots+\alpha_n\ge 0$, $\beta_1+\cdots+\beta_n\ge 0$.
 Let us assume $\sq{\x}\sim\sq{\y}$ in $t\LL$. Then, $ t\a+(\alpha_1,\ldots,\alpha_n)^\top\equiv t\b+(\beta_1,\ldots,\beta_n)^\top \;\mmod{t\M}$, that is, $t(\a-\b-\M\llambda)=(\beta_1-\alpha_1,\ldots,\beta_n-\alpha_n)^\top$ for some $\llambda\in\Z^n$. As $|\beta_i-\alpha_i|\in [0,t)$ for any $i=1,\ldots,n$, such an identity cannot hold except for the case $\beta_1-\alpha_1=\cdots=\beta_n-\alpha_n=0$. But, if so, we have either $\sq{\a}\sim\sq{\b}$ in $\LL$ with $\a\neq \b$, which contradicts Definition~\ref{def:hyperL} for $\LL$ to be a hyper-L, or $\a=\b$ and $\x=\y$ contradicting the first assumption.
\end{proof}
%
%

\begin{theorem}[The Dilating Method]
\label{teo:dilmet}
For an integer $t\geq1$,
\begin{itemize}
\item[$(a)$]
$\LL$ is a hyper-L of $G_{\M}$ if and only if $t\LL$ is a hyper-L of $G_{t\M}$.
\item[$(b)$]
$k(t\LL)=t(k(\LL)+n)-n$.
\item[$(c)$]
$\HH$ is an MDD (minimum distance diagram) of $G_{\M}$ if and only if $t\HH$ is an MDD of $G_{t\M}$.
\item[$(d)$]
$k(G_{t\M})=t(k(G_{\M})+n)-n$.
\item[$(e)$]
When applying the Dilating Method the set of generators is preserved.
\end{itemize}
\end{theorem}
\begin{proof}
$(a)$ The set $t\LL=\{t\sq{\a}:\sq{\a}\in\LL\}$ fulfills property $(ii)$ of Definition~\ref{def:hyperL} provided that $\LL$ is a hyper-L. Now we have to check property $(i)$. Notice that, for the digraph $G_\M$, the map $\phi$ is the `identity' map $\phi(\x)=x_1\e_1+\cdots+x_n\e_n\in\Z^n/\M\Z^n$ for $\x\in\N^n$.
By Lemma~\ref{lem:grid} we have $\sq{\x}\not\sim\sq{\y}$ for any pair of unit cubes $\sq{\x},\sq{\y}\in t\LL$. Assume $N=vol(\LL)=|\det\M|$. Then, $vol(t\LL)=t^nN=|\det(t\M)|$ and condition $(i)$ holds.
Analogous arguments can be used to show that $\LL$ is a hyper-L of $G_{\M}$ whenever $t\LL$ is a hyper-L of $G_{t\M}$.

The generators of $G_{t\M}$ are the same as those of $G_{\M}$. This fact is remarked in the proof of $(e)$.

$(b)$ Let us assume $k(\LL)=\norm{\a}$ for some $\sq{\a}\in\LL$. Then, by construction of $t\LL$, we have $k(t\LL)=\norm{t\a+(t-1,\ldots,t-1)^\top}=\sum_{i=1}^n(ta_i+t-1)=t\norm{\a}+n(t-1)=t(k(\LL)+n)-n$.

$(c)$ Assume $\HH$ is an MDD of $G_{\M}$. From $(a)$, the set of cubes $t\HH$ is a hyper-L of $G_{t\M}$. Now we have to check the optimal property of Definition~\ref{def:MDD},
\[
\norm{\a}=\min\{\norm{\x}:\x\in\N^n, \x\equiv\a\;\mmod{t\M}\}\quad\textrm{ for each }\sq{\a}\in t\HH.
\]
There are two kinds of unit cubes in $t\HH$. Namely,
\begin{itemize}
\item[$(c.1)$] $\sq{t\a}$ with $\sq{\a}\in\HH$;
\item[$(c.2)$] $\sq{t\a+(\alpha_1,\ldots,\alpha_n)^\top}$ with $\sq{\a}\in\HH$, $0\leq\alpha_1,\ldots,\alpha_n<t$, and $\alpha_1+\cdots+\alpha_n>0$.
\end{itemize}
Let us fix the value $t$ and consider a type-$(c.1)$ cube $\sq{t\a}\in t\HH$. Those $\x\in\N^n$ equivalent to $t\a$ in $\Z^n/t\M\Z^n$ are $\x=t\a+t\M\llambda$ for $\llambda\in\Z^n$. Then,
\begin{align*}
\min\{\norm{\x}:\x\in\N^n, \x\equiv t\a\;\mmod{t\M}\}&=\min\{\norm{t\y}:\y\in\N^n, \y\equiv\a\;\mmod{\M}\}\\
&=t\min\{\norm{\y}:\y\in\N^n, \y\equiv\a\;\mmod{\M}\}\\
&=t\norm{\a}=\norm{t\a},
\end{align*}
where the last line is a consequence of $\HH$ to be a minimum distance diagram and $\sq{\a}\in\HH$.

Consider now a type-$(c.2)$ cube $\sq{t\a+(\alpha_1,\ldots,\alpha_n)^\top}\in t\HH$. Set $\b_t=t\a+\bm{\alpha}$ with $\bm{\alpha}=(\alpha_1,\ldots,\alpha_n)^\top$. Assume that there is some $\y_t\equiv\b_t\;\mmod{t\M}$ with $\y_t\in\N^n$ and $\norm{\y_t}<\norm{\b_t}$ (and $\y_t\notin t\HH$). Then, $\y_t=\b_t+t\M\llambda$ for some $\llambda\in\Z^n$ and $\y_t=t(\a+M\llambda)+\bm{\alpha}=\c_t+\bm{\alpha}$ where $\c_t$ must have all its entries positive for being $\y_t\in\N^n$ (all the entries of $\c_t$ are multiple of $t$ and $0\leq\alpha_i<t$ for all $i$) and $\c_t\equiv t\a\;\mmod{t\M}$. Then,
\[
\norm{\y_t}=\norm{\c_t+\bm{\alpha}}=\norm{\c_t}+\norm{\bm{\alpha}}\geq\norm{t\a}+\norm{\bm{\alpha}}=\norm{t\a+\bm{\alpha}}=\norm{\b_t},
\]
which is a contradiction. Therefore, $t\HH$ is an MDD of $G_{t\M}$.

Assume now $t\HH$ is an MDD of $G_{t\M}$. Thus, $t\HH$ is a hyper-L of $G_{t\M}$ and, by  $(a)$, $\HH$ is a hyper-L of $G_{\M}$. Moreover, from $\sq{t\a}\in t\HH$ for each $\a\in\HH$, we have
\[
\norm{t\a}=\min\{\norm{\x}:~\x\in\N, \x\equiv t\a\;\mmod{t\M}\}=\min\{\norm{t\y}:~\y\in\N, \y\equiv\a\;\mmod{\M}\}
\]
that is equivalent to $\norm{\a}=\min\{\norm{\y}:~\y\in\N, \y\equiv\a\;\mmod{\M}\}$. So $\HH$ is also an MDD of $G_{\M}$.

$(d)$ Finally, the equality $k(G_{t\M})=t(k(G_{\M})+n)-n$ is a direct consequence of $(b)$ and $(c)$.

$(e)$
\label{lem:utmv}
Consider $\M\in\Z^{n\times n}$ with Smith normal form $\S=\U\M\V$. Then, as a direct consequence of the matrix product, the Smith normal form of $t\M$ is
\begin{equation}
t\S=\U(t\M)\V\label{eq:smtS}
\end{equation}
for an integral value $t\geq1$, which implies the result.

Moreover, following \eqref{eq:isom}, we can assume that $G_{\M}\cong\Cay(\Z^n/\M\Z^n,\{\u_1,\ldots,\u_n\})$ with $\u_1,\ldots,\u_n$ being the column vectors of $\S$. Then, using equality \eqref{eq:smtS}, it follows that $G_{t\M}\cong\Cay(\Z^n/t\M\Z^n,\{\u_1,\allowbreak\ldots,\u_n\})$ and so, the set of generators is preserved for any $t\geq1$.

\end{proof}

This result allows us to obtain infinite families of Cayley digraphs with the same density.

\vspace*{3mm}

\begin{corollary}
Let  $G_1=G_{\M}$ be an initial Cayley digraph of order $N$, degree $n$ and diameter $k$  Then,  the Cayley digraph $G_t=G_{t\M}$ has the same density as $G_1$ for any  integral value $t\geq1$.
\end{corollary}
\begin{proof}
 From Theorem~\ref{teo:dilmet}, we know that if $G_1$ has  minimum distance diagram $\HH$. then the dilates $t\HH$ are minimum distance diagrams related to the digraphs, of order $N_t=t^nN$ and diameter $k(G_{t\M})=t(k+n)-n$. Thus, their densities satisfy
\begin{equation}
\delta(G_{t\M})=\frac{N_t}{(k(G_{t\M})+n)^n}=\frac{t^nN}{(t(k(G_{\M})+n))^n}=\frac{N}{(k(G_{\M})+n)^n}=\delta(G_{\M}).\label{eq:density}
\end{equation}
\end{proof}

Thus,  the Dilating Method for obtaining dense families of Abelian Cayley digraphs results in the following procedure:
\begin{enumerate}
\item[{\bf 1.}]
Choose an initial (dense) digraph $G_1$ (perhaps found by computer search).
\item[{\bf 2.}]
Find any related minimum distance diagram $\HH$ of $G_1$.
\item[{\bf 3.}]
Use the fact that $\HH$ tessellates $\R^n$ to find
the set of column vectors of $\M$, $C_{\M}=\{\m_1,\ldots,\allowbreak\m_n\} \subset \Z^n$, such that $G_1\cong G_{\M}$.
\item[{\bf 4.}]
Apply the Dilating Method to obtain an infinite family of (dense) digraphs $G_t=G_{t\M}$.
\end{enumerate}
Notice that the steps {\bf 1},{\bf 2}, and {\bf 3} are done by using computer search. We remark again that all the elements of the infinite family of digraphs $\{G_t\}_{t\geq1}$ have the same density $\delta(G_1)$.
Note also that the Dilating Method generates an infinite family of non-cyclic Abelian Cayley digraphs excepting, perhaps, the initial one which depends on the selected digraph to apply the method.

\section{New dense families}
The criterion for a Cayley digraph $G_\M$ to be dense is not established and it is applied in the sense that $\delta(G_{\M})$ is as large as possible. This criterion is closely related to the degree-diameter problem for these digraphs. The parameter $\alpha=\alpha(G_{\M})$, where
\begin{equation}
\label{alpha}
N(G_{\M})=\alpha(G_{\M}) k(G_{\M})^d+O(k(G_{\M})^{d-1}),
\end{equation}
is also taken into account for a given infinite family of digraphs of degree $d$ and diameter $k(G_\M)=k$. In this context, notice that by \eqref{density} $\alpha(G)=\lim_{k\rightarrow \infty} \delta(G)$.

In fact, for a fixed degree $d=n$, several authors have proposed some infinite families of digraphs with good related values of $\alpha$, see Table~\ref{tab:d3}.  All these proposals are concerned with finite cyclic groups. The value $\alpha=0.0807$ applies only for diameters $k=22t+12$ with $t\not\equiv2,7\pmod{10}$. The case $\alpha=0.084$  applies only for diameters $k\equiv2\pmod{30}$.

\begin{table}[h]
\centering
\begin{tabular}{|l|l|l|}\hline
Paper&$\alpha$\\\hline\hline
G\'omez, Guti\'errez \& Ibeas  \cite{GoGuIb:07} (2007)&$0.037$\\
Hsu \& Jia  \cite{HsJi:94} (1994)&$0.062$\\
Aguil\'o, Fiol \& Garcia \cite{AgFiGa:97} (1997)&$0.074$\\
Chen \& GU  \cite{ChGu:92} (1992)&$0.078$\\
Aguil\'o  \cite{Ag:99} (1999)&$0.080$\\
Aguil\'o, Sim\'o \& Zaragoz\'a \cite{AgSiZa:01} (2001)&$0.084$\\\hline
\end{tabular}
\caption{Several proposals for degree $d=3$.}
\label{tab:d3}
\end{table}

Moreover, a result of Dougherty and Faber \cite[Cor.~8.2]{DoFa:04} proved the existence of Abelian Cayley digraphs of degree $d=3$ and order `at least' $N=0.084k^3+O(k^2)$ for all $k$ (and, hence, with $\alpha=0.084$). However, as far as we know, there is no explicit infinite family satisfying these conditions.


Finally, it is also worth mentioning the work of R\"odseth \cite{Rod96} on weighted loop networks, who gave sharp lower bounds for the diameter and mean distance for degree $d=2$ and general bounds for degree $d=3$.

The dense family given by Aguil\'o, Sim\'o, and  Zaragoz\'a \cite{AgSiZa:01} can be extended to a more general one. Using the same notation as in \cite{AgSiZa:01}, take the integral matrix $\M(m,n)$ given by
\begin{equation}
\M(m,n)=\left(\begin{array}{crr}
n&m&-2m-2n\\
3n+m&m&m+2n\\
2n&-m&m+n
\end{array}\right).\label{eq:mat-mn}
\end{equation}

\begin{proposition}[\cite{AgSiZa:01}]
\label{pro:old}
Consider the Cayley digraph $G_{m,n}=G_{\M(m,n)}$. Then, the order $N_{m,n}$ and diameter $k_{m,n}$ of $G_{m,n}$ are given by
\begin{align}
N_{m,n}&=m^3+12m^2n+14mn^2,\nonumber\\
k_{m,n}&\leq\max\{m+8n-3,3m+4n-3,5m-3\}.\label{eq:pro-old-k}
\end{align}
\end{proposition}

The case $\M(2,1)$ is given in Example~\ref{ex:M84}, with $N_{2,1}=84$ and $k_{2,1}\leq7$. The diameter $k=7$ is the minimum diameter a digraph $\Cay(\Z_{84},\{a,b,c\})$ can achieve. In \cite[Proposition~3]{AgSiZa:01} it is stated that, for $x\equiv0\;\mmod{3}$, the digraph $G_{\M(2x+1,x)}$ is isomorphic to a cyclic Cayley digraph and an explicit family is given. In the following result we extend this family to any value of $x$.

\begin{proposition}
\label{pro:bona1}
The Cayley digraph
\[
G_{\M(2x+1,x)}=\Cay(\Z_{84x^3 + 74x^2 + 18x + 1},\{-21x^2 - 15x - 2,21x^2+8x,-42x^2-23x-3\})
\]
has diameter $k_{2x+1,x}\leq10x+2$, for all integral value $x\geq1$.
\end{proposition}
\begin{proof}
The result follows from the Smith normal form decomposition of the matrix $\M(2x+1,x)$ in \eqref{eq:mat-mn}
\[
\S_x=\diag(1,1,84x^3 + 74x^2 + 18x + 1)=\U_x\M(2x+1,x)\V_x
\]
with unimodular matrices
\begin{align*}
\U_x&=
\left(\begin{array}{rcl}
-1 & 1 & -2 \\
-3 \, x - 1 & 3 \, x & -6 \, x - 1 \\
-21 \, x^{2} - 15 \, x - 2 & 21 \, x^{2} + 8 \, x & -42 \, x^{2} - 23 \, x - 3
\end{array}\right),\\
\V_x&=
\left(\begin{array}{rrr}
1 & 1 & -12 \, x^{2} - 10 \, x - 2 \\
0 & -1 & 12 \, x^{2} + 6 \, x + 1 \\
0 & 1 & -12 \, x^{2} - 6 \, x
\end{array}\right),
\end{align*}
and the upper bound \eqref{eq:pro-old-k} for the diameter in Proposition~\ref{pro:old}.
In fact, an alternative way to check that the first digraph$G_{\M(3,1)}=\Cay(\Z_{84},\{-21x^2 - 15x - 2,21x^2+8x,-42x^2-23x-3\})$
\end{proof}

Notice that this proposition does not provide an infinite family with density $0.084$.  Instead,
this value corresponds to the parameter $\alpha$, which, as commented above, is an asymptotic density when the diameter or the degree tend to infinity.
So, from the values of the order $N(x)=84x^3+74x^2+18x+1$ and diameter $k=k(x)=10x+2$, we get that $N(x)=\frac1{1000}(84k_{2x+1,x}^3+236 k(x)^2-152k_{2x+1,x}-312)=0.084 k(x)^3+O(k(x)^2)$ ($\alpha=0.084)$. Thus, if we look at the density of the elements of this family, we have
$\delta(x) = \frac{N(x)}{(k(x) + 3)^3}$.
Thus, there are no elements with density exactly $0.084$, as  $\delta(x)$ is strictly increasing from $\delta(1)=0.052444...$.

Now we give other dense families using the Dilating Method of the previous section.
More precisely, Proposition~\ref{pro:dil3} deals with a dense infinite family of digraphs of degree $d=3$ whereas Proposition~\ref{pro:diln} gives an infinite family for any degree $d\geq 2$.

Consider the integral matrices given in Example~\ref{ex:M84}, $\S$, $\U$, $\M$ and $\V$. The hyper-L $\HH$ given in Example~\ref{ex:L84} is a minimum distance diagram related to the Cayley digraph $G_{\M}$. This fact can be checked by computer. Now we can apply Theorem~\ref{teo:dilmet} to obtain the following result.
\\
\begin{proposition}
\label{pro:dil3}
The family of Abelian Cayley digraphs
\[
G_t=\Cay(\Z_t\oplus\Z_t\oplus\Z_{84t},\{(1,10,-38),(0,1,-3),(0,-2,7)\})
\]
has order $N_t=84t^3$ and diameter $k_t=10t-3$ for any integral value $t\geq1$.
\end{proposition}

Clearly, from (\ref{eq:density}), the density of this family is the constant value $\delta_t=0.084$ and, from $t=\frac{k_t+3}{10}$, in this case, we also obtain the parameter $\alpha_t=0.084$. It can be checked that the known Cayley digraphs of maximum density for degree $d=3$, listed in Section \ref{preli}, are contained in such a family. That is, $F_0\cong G_1$, $F_1\cong G_2$, $F_2\cong G_3$ and $F_3\cong G_4$.

Recall that, as commented above, Dougherty and  Faber \cite[Cor. 8.2]{DoFa:04} proved a existence (but no constructive) result about the existence of an infinite family of Abelian Cayley digraphs with degree $d=3$, order $N=0.084k^3+O(k^2)$, and density
$$
\delta(k)=0.084\left(\frac{k}{k+3}\right)^3+ O\left(\frac{1}{k}\right)
<0.084
$$
for all $k$.
In that paper, the authors derive their main results by using
the following notation:
$S_k=\{x\in \Z^d: \|x\|_1\le k\}$, and
$S_k'=S_k\cap P$, where $P$ is the positive octant of $\R^d$.
Then,  it is stated that $S_k'$ is a ``covering" of the lattice $\Z^d$, that is, $S_k'+N=\Z^d$, so that, in contrast with our method, this is not a proper `tiling' of $\Z^d$, as overlapping is allowed.
(Really, in \cite[Lem. 3.1]{DoFa:04}, it is said that we can find $T_k'\subset S_k'$ that tessellates $\Z^d$. However, Proposition 4.2(b) of the same paper deals with $S_k'$ and, hence, with possible overlapping.)

\begin{figure}[h]
\centering
\includegraphics[width=0.9\linewidth]{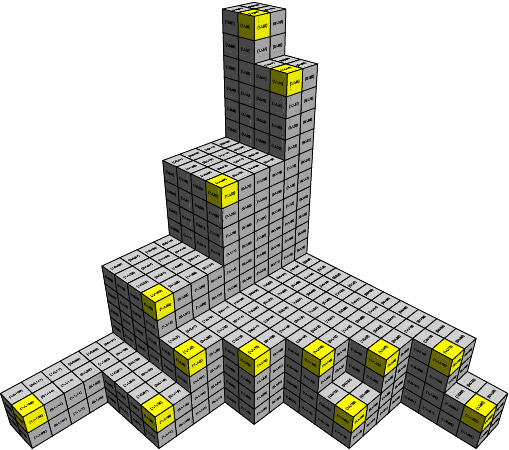}
\caption{Minimum distance diagram related to $G_{2}$ of Proposition~\ref{pro:dil3}.}
\label{fig:MDD-G2}
\end{figure}

\begin{example}
\label{ex:MDD-t2}
An MDD related to $G_2$ of Proposition~\ref{pro:dil3} is shown in Figure~\ref{fig:MDD-G2}. From the geometrical point of view, this diagram can be seen as the $2$-dilate of the diagram $\HH_1$ related to $G_1$ of Figure~\ref{fig:L84}, where each unitary cube has been dilated into four regular unitary cubes.
The 13 cubes with maximum norm $17=k(G_2)$ are $\sq{1,1,15}$, $\sq{1,3,13}$, $\sq{5,3,9}$, $\sq{9,3,5}$, $\sq{1,13,3}$, $\sq{3,11,3}$, $\sq{5,9,3}$, $\sq{7,7,3}$, $\sq{9,5,3}$, $\sq{1,15,1}$, $\sq{5,11,1}$, $\sq{11,7,1}$, and $\sq{15,1,1}$. Each cube $\sq{a_1,a_2,a_3}$ corresponds to the vertex
$$
\left(\begin{array}{rrr}1&0&0\\10&1&-2\\-38&-3&7\end{array}\right)
\left(\begin{array}{c}a_1\\a_2\\a_3\end{array}\right)\ \in\ \Z_2\oplus\Z_2\oplus\Z_{168}.
$$
\end{example}

In the following proposition we use the notation $\Z_m^{d-1}=\Z_m\oplus\stackrel{(d-1)}{\cdots}\oplus\Z_m$.

\begin{proposition}
\label{pro:diln}
Consider $A_d=\{(1,1,1,\ldots,1),(1,2,1,\ldots,1),\ldots,(1,1,\ldots,1,2)\}\subset\Z^d$. Then, the Cayley digraph $G_{d,t}=\Cay(\Z_{t}\oplus\Z_{t(d+1)}^{d-1},A_d)$ has diameter $k_{d,t}=t{d+1\choose 2}-d$, for $d\geq2$ and $t\geq1$.
\end{proposition}
\begin{proof}
Consider the family of digraphs given in \cite[Theorem~3.2]{AgFiPe:16a}, $F_{d}=\Cay(\Z_{d+1}^{d-1},B_d)$ with $B_d=\{(1,1,\ldots,1),(2,1,\ldots,1),\ldots,(1,\ldots,1,2)\}\subset\Z^{d-1}$ and $k(F_d)={d\choose2}$, defined  for $d\geq2$. Use the digraph isomorphism $F_d\cong G_{d,1}$ and apply the Dilating  Method to $G_{d,1}$ for $t\geq1$. From Theorem~\ref{teo:dilmet}$(d)$ it follows that $k(G_{d,t})=t(k(F_d)+d)-d=t{d+1\choose2}-d$.
\end{proof}

\begin{figure}[h]
\centering
\includegraphics[width=0.23\linewidth]{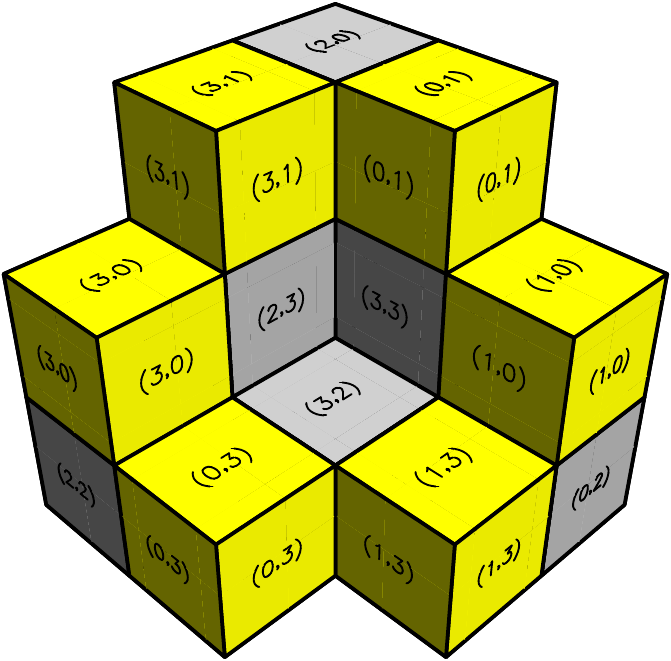}
\hspace{0.07\linewidth}
\includegraphics[width=0.5\linewidth]{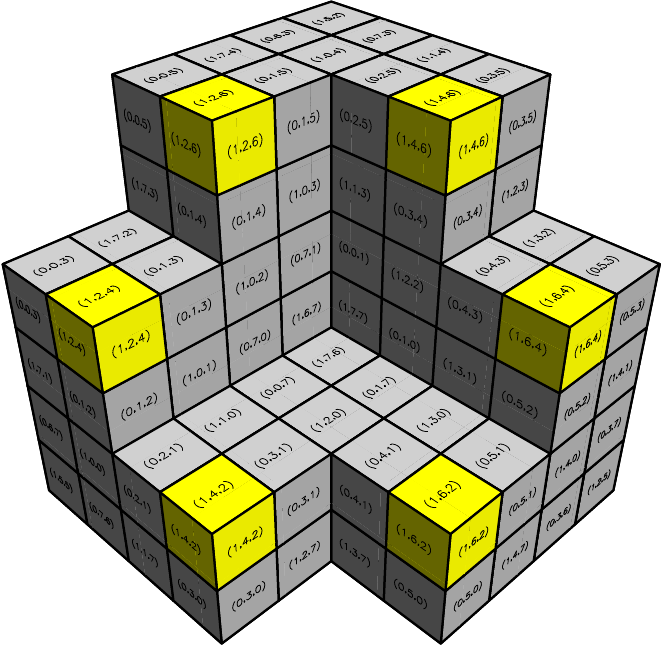}
\caption{Minimum distance diagrams related to $G_1=\Cay(\Z_4\oplus\Z_4,\{(1,1),(2,1),(1,2)\})$ and $G_2=\Cay(\Z_2\oplus\Z_8\oplus\Z_8,\{(1,1,1),(1,2,1),(1,1,2)\})$ of Proposition~\ref{pro:diln}.}
\label{fig:prop4}
\end{figure}

Notice the identity $\alpha(G_{d,t})=\delta(G_{d,t})=\frac1{d+1}(2/d)^d$, for all $t\geq1$.

By using the Stirling's formula, the lower bound expression $\lb(d,k)$ in (\ref{eq:df}),  gives
\[
\lb(d,k)\sim \frac{c}{\sqrt{2\pi}}e^{d-\frac{3}{2}\ln d-(\ln\ln d)(1+\log_2 e)}\left(\frac{k}{d} \right)^d+ O(k^{d-1}),
\]
with the multiplicative factor of $\left(\frac{k}{d} \right)^d$ being
\[
\frac{c}{\sqrt{2\pi}}e^{d-\frac{3}{2}\ln d-(\ln\ln d)(1+\log_2 e)}\sim \frac{c}{\sqrt{2\pi}}e^{d-\frac{3}{2}\ln d}.
\]
Notice that $N(G_{d,t})=t^d(d+1)^{d-1}$ and $k_{d,t}=k(G_{d,t})=t{{d+1}\choose2}-d$ hold for all $t$. Thus, as $t=\frac{k_{d,t}+d}{{{d+1}\choose2}}$, we get
\[
N(G_{d.t})=\frac{2^d}{d+1}\left(\frac{k_{d,t}}d+1\right)^d=\frac{2^d}{d+1}\left(\frac{k_{d,t}}{d} \right)^d+ O(k_{d,t}^{d-1}),
\]
with the multiplicative factor of $\left(\frac{k_{d,t}}{d} \right)^d$ being $e^{d\ln2-\ln(d+1)}$. To be compared with \eqref{eq:df}.

\section{Final comments}
Although the diameter of a digraph is a good measure of its metric properties, it is only an `absolute' measure. That is, two digraphs having the same order can be compared in terms of the best metric properties through the diameter. Roughly speaking, the lower the diameter, the better. However, when we want to compare two graphs with the same (optimal) diameter, but different orders, we need a ratio parameter   such as the density considered in this work.

Moreover, by using the density, we believe that a deeper knowledge of the tight lower bound expression, $\lb(d,n)$ in \eqref{eq:df}, can be reached. In this sense, when $d$ is fixed, computer search would complement and give light to this future work.

\end{document}